\documentclass[10pt,a4paper]{article}
\usepackage{cite}
\usepackage[T1]{fontenc}
\usepackage[utf8]{inputenc}
\usepackage{graphicx}
\usepackage{amsmath}
\usepackage{mathtools}
\usepackage{natbib}
\usepackage{authblk}

\usepackage{amsmath,amssymb}
\usepackage{paralist,enumitem}
\usepackage{verbatim}
\usepackage[right]{eurosym} 
\usepackage{url}
\usepackage{abstract}

\usepackage{algorithm2e}
\RestyleAlgo{ruled}

\usepackage{tikz}
\usetikzlibrary{decorations,decorations.pathmorphing,decorations.markings,decorations.text,%
decorations.footprints,decorations.fractals,decorations.pathreplacing,decorations.shapes,backgrounds,fit,%
trees,shapes.geometric,calc,shapes.symbols,patterns,positioning,shadows}
\tikzset{vertex/.style={circle,draw=blue!50,fill=blue!15,thick,minimum size=5mm,inner sep=0.5mm}} 

\tikzset{block/.style={vertex,rectangle}}
\tikzset{squigg/.style=decorate,decoration={snake,amplitude=0.5mm,segment length=1.5mm}} 
\tikzset{dreieck/.style={shape=isosceles triangle,%
  isosceles triangle apex angle=60,shape border rotate=90,fill=black!15,draw=black,thick,minimum width=1cm}}
\tikzset{>=stealth}
\usepackage{xcolor}

\newcommand{\R} {\mathbb{R}}

\begin{document}

\title{Modeling Minimum Cost Network Flows With Port-Hamiltonian Systems}
\date{}

\author{Onur Tanil Doganay\thanks{Corresponding author} }
\author{Kathrin Klamroth}
\author{Bruno Lang}
\author{Michael Stiglmayr}
\author{Claudia Totzeck}
\affil{University of Wuppertal, School of Mathematics and Natural Sciences, IMACM, Gaußstr.~20, 42119 Wuppertal, Germany\\
\{doganay,klamroth,lang,stiglmayr,totzeck\}@uni-wuppertal.de
}

\maketitle  

\begin{abstract}
We give a short overview of advantages and drawbacks of the classical formulation of minimum cost network flow problems and solution techniques, to motivate a reformulation of classical static minimum cost network flow problems as optimal control problems constrained by port-Hamiltonian systems (pHS).  
The first-order optimality system for the port-Hamiltonian system-constrained optimal control problem is formally derived. Then we propose a gradient-based algorithm to find optimal controls.
The port-Hamiltonian system formulation naturally conserves flow and supports a wide array of further modeling options as, for example, node reservoirs, flow dependent costs, leaking pipes (dissipation) and coupled sub-networks (ports). They thus provide a versatile alternative to state-of-the art approaches towards dynamic network flow problems, which are often based on computationally costly time-expanded networks. We argue that this opens the door for a plethora of modeling options and solution approaches for network flow problems. 
\end{abstract}

%%%%%%%%%%%%%%%%%%%%%%%%%%%%%%%%%%%%%%%%%%%%%%
\section{Introduction}

Minimum cost flow problems have an array of applications in operations research and thus have been extensively investigated since the 1960s \cite{ford62flows,ahuja93network,cruz22survey,chen22maximum}. While the original formulation considers a static situation and constant flow costs, there are extensions to the dynamic case with flow dependent costs. 
However, linear network flow models require extensive reformulation and linearization techniques to approximate time dynamics and/or non-linearity of flow costs. To avoid these reformulations which generally lead to a dramatic increase of the network size (in terms of the number of nodes and edges) we suggest a port-Hamiltonian (pHS) formulation of network flows which intrinsically covers non-linear costs and time-dynamic flows.

The pHS reformulation is motivated by \cite{schaft13port}. There network dynamics are modeled by ordinary differential equations. In particular, the dynamics of the flow along the edges of the network and dynamics of reservoirs at each node of the network are described. The approach is therefore versatile and allows for generalizations far beyond the static network flow case. In the following we discuss the simplest case of static dynamics and motivate further directions of research. Due to the general structure of ordinary differential equations we employ well-known techniques from optimization with differential equations leading us to an adjoint-based gradient descent method.

The remainder of this article is organized as follows. In Section~\ref{sec:mcnf} we summarize classical approaches to static linear minimum cost network flow problems \eqref{eq:MCNF} and some extensions to non-linear time-dynamic situations. In Section~\ref{sec:phs} we present a port-Hamiltonian system-model for network flows, which can be adapted easily to allow for non-linear cost functions and reservoirs in the nodes. Moreover, we formally derive an adjoint-based projected gradient method to compute  solutions for the minimum cost flow problem in port-Hamiltonian system formulation. Section~\ref{sec:concl} concludes the paper with some remarks to ongoing work and future perspectives.

\section{Classical formulation of minimum cost network flows and solution techniques}\label{sec:mcnf}

The task of minimum cost network flow problems is to find a cheapest transportation of flow through a predefined network that satisfies given supply, demand and capacity constraints. A typical application is the delivery of articles purchased online to customers. The nodes of the graph modeling  warehouses are assigned with the supply of goods and the demand nodes model customers purchasing a certain amount of the goods. 

In the simplest problem formulation, it is assumed that the goods are shipped directly from the warehouses to the customers leading us to so-called transportation problems. Capacity constraints limit the amount of goods that can be shipped along the routes. The more restrictive the capacity constraints, the more edges will be involved in the optimal flow. 

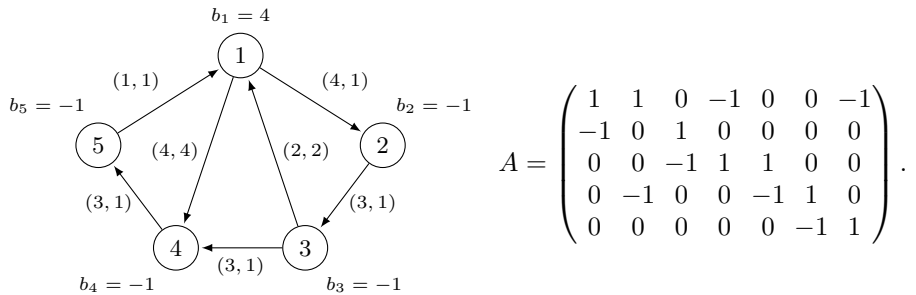
\begin{figure}[htb]
\begin{minipage}{0.54\textwidth}
\centering
      \begin{tikzpicture}[scale=0.85,-latex,>=stealth,shorten >=1pt,auto,node distance=2cm, main node/.style={circle,draw,font=\small\normalfont},
	every label/.style={font=\footnotesize}]

	\node[main node, label=above:{\scriptsize\(b_1=4\)}] (1) at (0,3) {$1$};
	\node[main node, label={80:{\scriptsize\(b_2=-1\)}}] (2) at (2.2,1.6) {$2$};
	\node[main node, label=-60:{\scriptsize\(b_3=-1\)}] (3) at (1,0) {$3$};
	\node[main node, label=-120:{\scriptsize\(b_4=-1\)}] (4) at (-1,0) {$4$};
	\node[main node, label=100:{\scriptsize\(b_5=-1\)}] (5) at (-2.2,1.6) {$5$};

	\path[every node/.style={font=\scriptsize}]
	  (1) edge node [above right] {\((4,1)\)} (2)
	      edge node [left] {\((4,4)\)} (4)
	  (2) edge node [pos=0.6,right] {\((3,1)\)} (3)
	  (3) edge node [right] {\((2,2)\)} (1)
	      edge node [below] {\((3,1)\)} (4)
	  (4) edge node [pos=0.4,left] {\((3,1)\)} (5)
	  (5) edge node [above left] {\((1,1)\)} (1);
      \end{tikzpicture}
\end{minipage}
\begin{minipage}{0.4\textwidth}
\[\setlength{\arraycolsep}{0.4ex}
A = \begin{pmatrix} 1 & 1 & 0 & -1 & 0 & 0 & -1 \\
         -1 &  0 & 1 & 0 & 0 & 0 & 0 \\
         0 & 0 & -1 & 1 & 1 & 0 & 0 \\
         0 & -1 & 0 & 0 & -1 & 1 & 0 \\
         0 & 0 & 0 & 0 & 0 & -1 & 1 \end{pmatrix}.
\]
\end{minipage}
\caption{The values next to the vertices denote the supply (positive) and the demand (negative) of each vertex. We assume that total supply equals total demand. The tuples $(u_e,c_e)$ at the edges denote their capacities and costs, respectively. An initial feasible flow is given by $x = (4,2,3,2,0,1,0)^\top$.}
\label{fig:uebung}
\end{figure}

When we track our parcels online, we see that they are sometimes reloaded at redistribution hubs, which are spread all over the country. These hubs can be modeled using transshipment nodes, and lead to more sophisticated problems but at the same time more efficient routing. See Figure~\ref{fig:uebung} for an illustrative example with one warehouse in node $1$ and four customers (nodes 2--5) each ordering one unit. 
Classical problem formulations and implementation of the solution methods often assume integer values of supply, demand and capacities. Let us recall the classical formulation of minimum cost network flow problems and introduce the notation that will be generalized in Section~\ref{sec:phs} to time-dependent dynamics.

\medskip
We consider a directed graph $G=(V,E)$ with $N_e$ edges and $N_v$ nodes. The connectivity of the graph is encoded in the incidence matrix $A \in \{-1,0,1\}^{N_v \times N_e}.$ Every column of \(A\) represents one edge of \(G\) and has exactly two non-zero entries, where $a_{ij}=1$ if edge \(j\) is an outgoing edge of node \(i\) and $a_{ij}=-1$ if \(j\) is an incoming edge of node $i$. The flow on edge $e=(i,j)\in N_e$ is denoted by $x_e(t)=x_{i,j}(t)\in\R$ and stored in a vector of flow values $x(t) = (x_e(t))_e$. Storing the supply/demand values in a vector \(b\in\R^{N_v}\) (supplies positive, demands negative and \(\sum_{v\in V}b_v=0\)), the flow conservation constraints can be easily formulated as \(A\,x=b\).
We assume thereby that the edges have upper and lower capacity limits \(x_l\) and \(x_u\), respectively, such that \(x_l\leq x \leq x_u\). Moreover, the transportation of one unit of flow on edge $e$ is associated with costs, \(c_e \ge 0\). Then, the minimum cost network flow problem is given by the following linear program

\begin{equation}\tag{MCNF}\label{eq:MCNF}
  \begin{array}{rr@{\extracolsep{0.75ex}}c@{\extracolsep{0.75ex}}l}
   \min &\multicolumn{3}{l}{\displaystyle \sum_{e \in E} c_e\, x_e}\\
   \text{s.\,t.} & A\,x &=& b\\
   & x_\ell &\leq &x \leq x_u,
  \end{array}
\end{equation}
where the first constraint accounts for the conservation property in each node including supply and demand nodes, with supplies and demands given by the vector $b\in\R^{N_v}$ with $\sum_{v\in V}b_v=0$. 

In general, any solution method for linear programs can be employed to solve \eqref{eq:MCNF}, but it has been shown that it pays off to exploit the additional structure given by the network. Indeed, the 
network simplex algorithm uses the residual network and the specific structure of the incidence matrix in each of its iterations and has become a well-established solution technique over the years. Its efficiency is due to the fact that both pivot selection and basis exchange can be implemented extremely efficiently. Since the basic variables of \eqref{eq:MCNF} form a spanning tree, the corresponding columns of the full-rank incidence matrix can be rearranged (by a reordering of rows/nodes and columns/edges) to form an upper triangular matrix. Moreover, the update of the flow values in each iteration of the network simplex is limited to edges of an (undirected) cycle which can be easily determined from the tree of basic edges. For a self-contained introduction to the network simplex and efficient implementation techniques see, e.g., \cite{ahuja93network,korte18combinatorial,cruz22survey,chen22maximum}. 

In contrast, time-dynamic minimum cost flow problems and problems with non-linear flow costs are much more difficult to solve. This includes generalizations of \eqref{eq:MCNF} with time- and/or flow-dependent costs, time-dependent capacities or time-dependent supplies and demands. To preserve the linear structure of \eqref{eq:MCNF} often reformulation techniques are applied that approximate these generalizations by a static \eqref{eq:MCNF} on an extension of the network structure with constant cost coefficients.

Time-dynamic network flows (as for example in the quickest flow problem \cite{fleischer07quickest}) are often modeled using so-called \emph{time-expanded networks}. Based on a time-discretization \(\Delta=(t_0, t_0+\delta, t_0+2\,\delta, \ldots )\), a copy of the set of nodes is considered for each time-step in a layer representing this time-step. Then nodes \(u^i\) and \(v^j\) in layers \(i\) and \(j\) are connected by an edge, if and only if there is an edge \(e=(u,v)\) in the original graph and a unit of flow needs time of \(\delta\,(j-i)\) to traverse edge \(e\). Using time-expanded networks dynamic network flow problems can be solved approximatively. However, this reformulation increases the network size (number of nodes and edges) by orders of magnitude. In particular, finer resolutions of the time-discretization and larger networks lead to prohibitively large instance sizes. For a general introduction to dynamic network flow problems, their modeling and solution methods see, e.g., \cite{pyakurel20network,kotnyek03annotated,pageni21survey}.

Network flow problems with non-linear flow dependent edge costs are often approximated by replacing each edge by multiple parallel edges with limited capacities and increasing cost coefficients. In the optimal solution the flow is first routed over the cheapest of the parallel edges and only if the flow value exceeds its capacity limit a more expensive parallel edge is used \cite{koehler05flows}. In the following we suggest an alternative approach to (dynamic) minimum cost network flow problems which can be easily generalized to  flow-dependent costs and time-dependent settings.

\section{Modeling minimum cost network flow as port Hamiltonian system}\label{sec:phs}
The static graph structure with flow conservation along the nodes can be naturally generalized to time-dependent dynamics using the port-Hamiltonian framework \cite{schaft13port}. In addition to the network and flows introduced above, we consider reservoirs in each node of the network. Analogous to the definition of the flows, we define potentials (or pressure) in the reservoirs by $\rho = (\rho_v)_v$. The interplay of flow and potentials can be modeled as a system of ordinary differential equations in form of a port-Hamiltonian system \cite{schaft13port}. In more detail, the dynamics read
\begin{alignat*}{3}
	&\frac{d}{dt} \rho = A \, x + u, \qquad \quad&\rho(0) = \hat \rho, \\
	&\frac{d}{dt} x = -A^\top \rho, &x(0) = \hat x.
\end{alignat*} 
The first equation models the change of the potential in the reservoirs, which is defined by the amount flowing into or out of the pipes, respectively. The direction is coded in the signs of the incidence matrix $A$. The second equation models the change of flow along the edges which is varying due to the potentials (or pressure) in the reservoirs. We emphasize that the model is flow conserving by construction. 
To see the port-Hamiltonian structure we define $z = (\rho,x)$ and rewrite the system as
\begin{equation}\tag{PHS}\label{eq:phs}
	\frac{d}{dt} z = \mathcal J(z) \, \nabla \mathcal H(z) + B \, u,  \qquad z(0) = \hat z, 
\end{equation} 
where 
\begin{equation}	    
\mathcal J(z) = \mathcal J = \begin{pmatrix} 0 & A \\ -A^\top & 0  \end{pmatrix}
\end{equation} 
encodes the connectivity of edges and nodes and
\[
\mathcal H(x,\rho) = \frac12 x^2 + \frac12\rho^2
\]
is the Hamiltonian of the dynamics and $B = \begin{pmatrix}I \\ 0 \end{pmatrix}$ defines the input action. For completeness, we note that the output of the system is given by $y = B^\top z = \rho.$ It does not play a role in the following, but will be considered in future work. 

Note that $\mathcal J$ is skew-symmetric, hence the conservation property of \eqref{eq:phs} can be shown by computing the time derivative of $\mathcal H$
\begin{equation}\label{eq:dissipativity}
    \frac{d}{d t} \mathcal H(x,\rho) = \big( \nabla \mathcal H, \frac{d}{d t} z \big) = \big( \nabla \mathcal H, \mathcal J(z) \, \nabla \mathcal H(z) + B \, u \big) = \big( \nabla \mathcal H, B \, u \big) = \big( y, u\big)
\end{equation}
In particular, this shows that for $u=0$ or $y=0$ the dynamics is conservative.
 
 We emphasize that the flow along the edges and the potential (or pressure) in the reservoirs of the port-Hamiltonian system are  time-dependent quantities, which allows for richer applications compared to \eqref{eq:MCNF}. Other degrees of freedom are the structure of the Hamiltonian $\mathcal H$ and the input as well as generalizations to state dependent incidence matrices $A(z).$
 
 \subsection{Relationship of port-Hamiltonian system dynamics and constraints of the minimum cost network flow problem}
 The first step to linking the constraints of the minimum cost network flow problem with the port-Hamiltonian system dynamics is to choose the initial conditions consistently. In fact, $\hat x$ is a feasible solution of \eqref{eq:MCNF}, if $A\hat x = b$ and $x_\ell \leq \hat x \leq x_u$. As the minimum cost network flow problem admits no reservoirs, we choose $\hat\rho \equiv 0$. To preserve the feasibility of $x$ over time, we choose the feedback $u = -Ax$. Altogether, this leads to a stationary solution $\rho \equiv 0$ and $x \equiv \hat x.$ Hence, $x(t)$ is a feasible flow of \eqref{eq:MCNF} for all $t \in [0,T]$. Moreover, the output satisfies $y\equiv 0$. Hence by the information of \eqref{eq:dissipativity} the dynamics is conservative, i.e., there is no dissipation and all flow is conserved. The details of the dynamics with feedback are given by
 \begin{subequations}\label{eq:constraint}
 \begin{align}
	\frac{d}{dt} z &= \mathcal J  \, \nabla \mathcal H(z)
	+ \begin{pmatrix}I \\ 0 \end{pmatrix} u,  \\ z(0) &= (0, \hat x)^\top,\\
	u&\coloneqq -Ax.
\end{align}
\end{subequations}
We emphasize that we denote by $\hat x$ the initial condition and by $x: [0,T]\rightarrow \R^{N_e}$ a feasible flow over time. Our goal is to identify $\hat x$ that is a feasible starting condition and yields, depending on the objective, a minimum cost flow $x$. 
In the following, we formulate an optimal control problem constrained by \eqref{eq:phs} that allows for generalizations of \eqref{eq:MCNF}.  The first aim is to provide a port-Hamiltonian system formulation of the static minimum cost network flow problem. In future work, the advantages of the port-Hamiltonian model are exploited to establish a framework for richer minimum cost flow problems.

\subsection{Minimum cost flow problem in port-Hamiltonian system formulation}
As we want to consider time-dependent problems in the long term, we aim for an adjoint-based optimization algorithm to solve the minimum cost flow problem in port-Hamiltonian system formulation. A straight-forward choice of the cost functional would be $\int_0^T c_e x_e(t) dt$ analogous to the static minimum cost flow problem. Note that this is inconvenient for our approach, as the linearization is independent of the state $z$ and therefore, no state information would enter the adjoint system. To overcome this issue, we propose the cost functional 
\begin{equation*}
	\min\limits_{\hat x}  \int_0^T \left(\sum_{e \in E} c_e\, x_e(t)\right)^2  dt.
\end{equation*}
Note that the linearization in arbitrary direction $h$ is now given by 
\begin{equation*}
    d_x \left( \int_0^T \left(\sum_{e \in E} c_e\, x_e(t)\right)^2  dt \right) [h] = \int_0^T 2 \left( \sum_{e \in E} c_e\, x_e(t) \right) c_e h_e(t) dt.
\end{equation*}
Now, the adjoint state depends on the current state of the flow as the prefactor is state dependent. However, as this is just a scaling of the cost vector,  the influence of the state is only marginal.

By construction of the feedback, the flow is constant over time, i.e., $x(t) = \hat x$ for all $t\in[0,T]$ and the reservoirs remain empty. Hence, 
\begin{equation*}%\label{eq:cost}
	\min\limits_{\hat x}  \int_0^T \left(\sum_{e \in E} c_e\, x_e(t)\right)^2  dt = \min\limits_{\hat x}  T \left(\sum_{e \in E} c_e\, \hat x_e \right)^2 .
\end{equation*}
As $c_e \ge 0$ for all $e\in E$ the optimal initial flow $\hat x$ coincides with the minimizer of \eqref{eq:MCNF}. From now on we interpret the initial condition as control, therefore we introduce the notation $x(0) = \hat x =w$. 

Altogether, this motivates to consider an optimal control problem in the port-Hamiltonian system framework:
\begin{equation*}
	\min\limits_{w \in \mathcal W_\text{ad}}  \int_0^T \left(\sum_{e \in E} c_e\, x_e(t)\right)^2  dt \quad \text{subject to } \, \eqref{eq:constraint} 
\end{equation*}
with set of admissible controls given by $$\mathcal W_\text{ad} = \{ w \in \mathbb R^{N_e} \colon Aw = b, \, x_\ell \le w \le x_u\}.$$ 

Having again the long-term perspective of richer time-dependent applications in mind, we propose to use a gradient-based algorithm. Therefore, the derivation of the optimality condition of the optimal flow problem is based on a Lagrangian approach with Lagrange function $\mathcal L$ and state operator $\bar{e}$, respectively, given by
\[
\mathcal L(z,\xi, w) = \int_0^T C(z(t))  \, dt - \langle \xi, \bar{e}(z,w) \rangle, \qquad \bar{e}(z,w) = \begin{pmatrix} \frac{d}{dt} z - \mathcal J \, \nabla \mathcal H(z) \\ x(0) - w \end{pmatrix},
\]
where $C(z(t)) =  \left(\sum_{e \in E} c_e\, x_e(t)\right)^2$ and the control $w$ is the initial condition of the flow.
Using standard methods \cite{hinze09opt}, the first-order optimality condition can be computed by solving $d\mathcal L = 0.$ This yields the first-order optimality system consisting of 
\begin{align*}
    &d_\xi \mathcal L = 0 \quad &\Leftrightarrow \qquad &\frac{d}{dt} z = \mathcal J \, \nabla \mathcal H(z), \quad z(0) = (0, \hat x)^\top \\ & & & \text{(state equation)} \\
    &d_z \mathcal L = 0 \quad &\Leftrightarrow \qquad &\frac{d}{dt} \xi = \mathcal J^\top \xi - C'(z), \quad\xi(T) = 0, \\ & & &  \text{(adjoint equation)} \\
    &(d_w \mathcal L, w- w_\mathrm{opt}) \ge 0 \quad &\Leftrightarrow \qquad &\Big(-\xi(0), w-w_\mathrm{opt}\Big) \ge 0 \quad \text{for all } w \in \mathcal W_\text{ad}  \\ & & &  \text{(optimality condition)}
\end{align*}
which can be employed for a projected gradient descent method to compute the optimal initial condition $w_\mathrm{opt}.$ Here, the projection on $\mathcal W_\text{ad}$ ensures that the gradient update leads to feasible flows; in more detail we compute a gradient candidate $g$ by solving the linear program
\begin{gather}\tag{projection}\label{eq:projection}
\min\limits_{\hat g} (-\xi(0), \hat g) \qquad 
\text{subject to} \qquad A\, \hat g = 0, \quad x_\ell- \hat x \le \hat g \le x_u- \hat x.
\end{gather}
In the next section we show that the optimality condition coincides with the gradient of the reduced cost functional, which allows us to formulate a  gradient descent algorithm for the minimum cost network flow problem.

\subsection{Derivation of the gradient}

By standard ODE theory \cite{Teschl}, \eqref{eq:constraint} admits a unique solution. This allows us to define the control-to-state map $S(w) = z$ and moreover the reduced cost functional $f(w) := C(S(w)).$ The state constraint is implicitly satisfied in the reduced formulation, hence the gradient descent algorithm proposed in the following only varies the initial condition. First recall that $\bar{e}(z,w)=0$ and hence $\bar{e}_z(z,w) S'(w) + \bar{e}_w(z,w) = 0.$  This is used to obtain
\begin{align*}
    \langle f'(w), h_w \rangle &= \langle C'(z), S'(w) h_w \rangle = \langle S'(w)^* C'(z), h_w \rangle \\
    &= \langle-\bar{e}_w(S(w),w)^* \bar{e}_y(S(w),w)^{-*} C'(z),h_w \rangle \\
    &= \langle \bar{e}_w(S(w),w)^* \xi,h_w \rangle \\
    &= \langle -\xi(0), h_w \rangle
\end{align*}
for an arbitrary direction $h_w \in \mathcal W_\text{ad}$. Here, we use the following notation: for an operator $T,$ we denote its adjoint by $T^*$ and the inverse of the adjoint by $T^{-*}.$ Note that the last equation coincides with the optimality condition for $h_w = w-w_\mathrm{opt}.$
Now, Riesz representation theorem allows us to identify the gradient $\nabla f(w).$ Since our control is the initial condition $\hat x\in \mathbb R^{N_e},$ we obtain $\nabla f(w) = -\xi(0),$ which we project to $\mathcal W_\text{ad}$ via \eqref{eq:projection}. We denote this in the following using the projection operator $\mathcal P_{\mathcal{W}_\text{ad}}$. Then we can proceed with the gradient descent algorithm. 

\medskip

\subsection{Gradient descent algorithm}
Note that both, the network simplex and the gradient descent algorithm, need to be initialized with a feasible initial solution $w\in\mathcal W_\text{ad}.$ We propose to compute a feasible initial flow by solving a maximum flow problem on an associated graph with an additional source node $s$ and an additional terminal node $t$. For this purpose, all original supply nodes are connected from the newly introduced starting node \(s\) by an edge with capacity equal to their respective supply values. Analogously, the demand nodes are connected to a new terminal node \(t\) by an edge with capacity equal to their respective demands. In the resulting network, all nodes are transshipment nodes (i.e., all supplies and demands are set to zero). Then any maximum \(s\)-\(t\) flow corresponds to a feasible flow in the original network. Maximum flows can be efficiently computed with the \textit{Algorithm of Ford-Fulkerson} or \textit{Preflow-Push methods}, see, e.g., \cite{ahuja93network,korte18combinatorial}. 

We have now all ingredients for the gradient descent algorithm to compute minimum cost network flows at hand. In fact, we start with an initial guess and solve the state equation. The initial guess and the corresponding state solution are data for the adjoint equation, which we solve backward in time. The solution of the adjoint at time zero is then projected via $\mathcal P_{\mathcal{W}_\text{ad}}$ to obtain a feasible update for the next interation. A pseudocode for the minimum cost network flow problem with port-Hamiltonian system-constraint is given in Algorithm~\ref{alg:gradientdescent}. To find an appropriate step size we propose to use the Armijo rule, see for example \cite{hinze09opt}.

\begin{algorithm}
\caption{Gradient descent for minimum cost network flow problem}\label{alg:gradientdescent}
\KwData{feasible initial guess $w,$ \\ incidence matrix $A,$ \\ other algorithmic parameters\;}
\KwResult{optimal initial condition $w$}
\medskip
solve state problem to get  $S(w)$\;
solve adjoint problem to get $\xi$ for given $S(w)$ and $w$\;
identify the gradient and project to $\mathcal W_\text{ad}$ to obtain $\mathcal P_{\mathcal W_\text{ad}}(\nabla f(w))$\;
\While{$|\mathcal P_{\mathcal W_\text{ad}}(\nabla f(w))|\geq \epsilon_\text{stop}$}
{   
    choose appropriate step size $\lambda$ with Armijo rule\; 
    $w \gets w - \lambda\,\mathcal P_{\mathcal W_\text{ad}}(\nabla f(w))$\;
    solve state problem to get  $S(w)$\;
solve adjoint problem to get $\xi$ for given $S(w)$ and $w$\;
identify the gradient and project to $\mathcal W_\text{ad}$ to obtain $\mathcal P_{\mathcal W_\text{ad}}(\nabla f(w))$\;
}
\end{algorithm}

Preliminary numerical results indicate that the time horizon on which the state and the adjoint system are solved can be very short. This is due to the enforced stationarity of the problem. In our preliminary tests, most of the computational effort is needed to solve the projection on $\mathcal W_\text{ad}$. Moreover, it can be seen that the port-Hamiltonian system solver chooses combinations of dicycle flows in the residual network in order to reduce the cost as fast as possible, while the standard simplex methods iteratively reduces costs along only one dicycle in the residual network in every iteration. The combination of several dicycle flows chosen in the port-Hamiltonian system is expected to be particularly advantageous for large scale problems since it has the potential to significantly reduce the number of iterations needed until an optimal flow is detected. For static problems, however, this advantage can not compensate for the significantly higher cost for computing the projection.

\section{Conclusions and perspectives}\label{sec:concl}

We plan to investigate the relationship in terms of well-posed\-ness of the minimization problems and the adjoint-based solution approach. Although the port-Hamiltonian system structure allows for more modeling features, such as reservoirs and time dependence, we see analogies between the static residual network that is used in the simplex method for network flow problems \cite{ford62flows,ahuja93network} and the projection to the feasible set \cite{troeltzsch10optimal} in the optimal control problem that we plan to investigate in detail. To this end, we aim to exploit the optimality conditions of the two problems. 

Based on these findings, a natural next step is to generalize the method to problems with nonlinear cost $c_e(x_e)$ and flow dependent incidence matrices $A(x)$. Note that this type of generalizations are straight forward in the port-Hamiltonian system setting. In contrast, the computational effort to solve the minimum cost flow problem changes significantly, as we have already discussed in Section~\ref{sec:mcnf}. In fact, in nonlinear settings the port-Hamiltoanian systems approach may be advantageous in terms of computational cost.

Other extensions could be the consideration of time-dependent supply and demand functions $b(t)$ and time-dependent capacities $x_\ell(t)$ and $x_u(t)$. Here, the idea of finding a stationary flow is infeasible. On the one hand, this leads us to the challenge of finding an appropriate dynamical minimum cost flow formulation. On the other hand, we need to model the control action in the optimal port-Hamiltonian system flows appropriately.

Further, we aim to take advantage of the additional features that come with the port-Hamiltonian system formulation of dynamics on graphs which can be applied in the modeling of real-world applications such as evacuation, transportation or production networks. Reservoirs that are naturally contained in the port-Hamiltonian system dynamics can be used to model waiting queues in evacuation scenarios as well as production networks. To generalize the model even more, we can consider abstract incidence matrices $A(z).$ Other degrees of freedom are the Hamiltonian $\mathcal H(z).$ Moreover, we can include dissipativity with the help of a matrix $R$ (to model, e.g., leaking pipes), control actions $B(z)u$ and corresponding ports $y$ leading to the port-Hamiltonian system dynamics given by 
\begin{align*}
\frac{d}{dt} z &= (\mathcal J-R) \,\nabla \mathcal H(z) + B(z)\, u, \\ z(0) &= \hat z, \\
y &= B(z)^\top \nabla \mathcal H(z).
\end{align*}
In particular, we want to use ports and inputs for the decomposition tasks applicable to large-scale networks and to couple several sub-networks. One application would be the modeling of evacuation of buildings with several identical floors. The sub-networks are identical on each floor and coupled through stairs and elevators.

\bibliographystyle{abbrvnat} 
\bibliography{PHS_networks}

\end{document}